# ON THE SHINTANI ZETA FUNCTION FOR THE SPACE OF BINARY TRI-HERMITIAN FORMS

Akihiko Yukie[1]

Oklahoma State University

ABSTRACT. In this paper, we consider the most non-split parabolic $D_4$ type prehomogeneous vector space. The vector space is an analogue of the space of Hermitian forms. We determine the principal part of the zeta function.

## Introduction

Throughout this paper, $k$ is a number field. Let $k_1$ be a cubic extension of $k$, and $k_2$ the Galois closure of $k_1$. Then either $k_2 = k_1$ is a cyclic cubic extension of $k$ or $k_2$ is an $\mathfrak{S}_3$–extension of $k$. In this paper, we consider the zeta function defined for a prehomogeneous vector space $(G, V)$ which becomes the $D_4$ case in [6] after extending the base field from $k$ to $k_2$. We consider the most non-split case here and the vector space is an analogue of the space of Hermitian forms except it depends on three two dimensional vectors. So we call this space the space of binary tri-Hermitian forms. We discussed the meaning of this case in [1].

The purpose of this paper is to determine the principal part of the zeta function. Our main result is Theorem (3.6). In §1, we define the prehomogeneous vector space in this paper, and prove a stratification of this vector space. Since our group is non-split, we cannot directly apply Kempf's theorem on the rationality of the equivariant Morse stratification [2], even though it does not seem so difficult to make necessary adjustments. However, since our representation is rather small, we will handle the stratification explicitly. We will probably consider such adjustments for general non-split prehomogeneous vector spaces in the future.

In §2, we define the zeta function, the Fourier transform, and the smoothed Eisenstein series. Our method for computing the principal part of the zeta function is similar to that of [3], [5], [7], [8]. In §3, we will analyze the terms which come from the set of unstable points, and prove a principal part formula.

For basic notion on adeles, see [4]. The ring of adeles, the group of ideles, and the adelic absolute value of $k$ are denoted by $\mathbb{A}$, $\mathbb{A}^\times$, $|\ |$ respectively. We fix a non-trivial character $<\ >$ of $\mathbb{A}/k$. The ring of adeles, the group of ideles, and the adelic absolute value of $k_1$ are denoted by $\mathbb{A}_{k_1}$, $\mathbb{A}_{k_1}^\times$, $|\ |_{k_1}$ respectively. Note that by the

[1]The author is partially supported by NSF grant DMS-9401391 and NSA grant MDA-904-93-H-3035

Typeset by $\mathcal{A}\mathcal{M}\mathcal{S}$-TEX

inclusion $\mathbb{A} \to \mathbb{A}_{k_1}$, an idele $(a_v)_v$ corresponds to the idele $(b_w)_w$ such that $b_w = a_v$ if $w$ is a place over $v$. Let $\mathbb{A}^1 = \{t \in \mathbb{A} \mid |t| = 1\}$, $\mathbb{A}^1_{k_1} = \{t \in \mathbb{A}_{k_1} \mid |t|_{k_1} = 1\}$. Identifying $k_1 \otimes \mathbb{A} \cong \mathbb{A}_{k_1}$ the norm map $N_{k_1/k}$ can be extended to a map from $\mathbb{A}_{k_1}$ to $\mathbb{A}$. It is known (see [4, p. 139]) that $|N_{k_1/k}(t)| = |t|_{k_1}$ for $t \in \mathbb{A}_{k_1}$. The set of positive real numbers is denoted by $\mathbb{R}_+$.

Suppose $[k : \mathbb{Q}] = n$. Then $[k_1 : \mathbb{Q}] = 3n$. For $\lambda \in \mathbb{R}_+$, $\underline{\lambda} \in \mathbb{A}^\times$ is the idele whose component at any infinite place is $\lambda^{\frac{1}{n}}$ and whose component at any finite place is 1. Also $\underline{\lambda}_{k_1} \in \mathbb{A}^\times_{k_1}$ is the idele whose component at any infinite place is $\lambda^{\frac{1}{3n}}$ and whose component at any finite place is 1. Clearly, $\underline{\lambda} = \underline{\lambda}^3_{k_1}$. Since $|\underline{\lambda}| = \lambda$ and $|\underline{\lambda}_{k_1}|_{k_1} = \lambda$, this means that $|\underline{\lambda}|_{k_1} = \lambda^3$.

If $X$ is a variety over $k$ and $R$ is a $k$–algebra, the set of $R$–rational points of $X$ is denoted by $X_R$. The space of Schwartz–Bruhat functions on $V_\mathbb{A}$ is denoted by $\mathscr{S}(V_\mathbb{A})$. If $f, g$ are functions on a set $X$, $f \ll g$ means that there is a constant $C$ such that $f(x) \le Cg(x)$ for all $x \in X$.

### §1 The space of binary tri-Hermitian forms

Let $G_1 = \mathrm{GL}(1)_k$, $G_2 = \mathrm{GL}(2)_{k_1}$, and $G = G_1 \times G_2$. We consider $G$ as a group over $k$. Let $\widetilde{V} = k_2^2 \otimes k_2^2 \otimes k_2^2$. We first fix a coordinate system for $\widetilde{V}$. Let $f_1 = \begin{pmatrix} 1 \\ 0 \end{pmatrix}$, $f_2 = \begin{pmatrix} 0 \\ 1 \end{pmatrix}$, and $e_{ijk} = f_i \otimes f_j \otimes f_k$ for $i, j, k = 1, 2$. Then $\{e_{ijk} \mid i, j, k = 1, 2\}$ is a basis for $\widetilde{V}$. So any element $x \in \widetilde{V}$ can be expressed as $x = \sum_{i,j,k} x_{ijk} e_{ijk}$, where $x_{ijk} \in k_2$ for all $i, j, k$.

We choose three different imbeddings $\sigma_1, \sigma_2, \sigma_3 : k_1 \to k_2$ over $k$. Let $H_1 = \mathrm{Gal}(k_2/k)$, and $H_2 = \mathrm{Gal}(k_2/k_1)$. Then $H_2$ is a subgroup of $H_1$ and $[H_1 : H_2] = 3$. So there exists a homomorphism $h : H_1 \to \mathfrak{S}_3$ such that $H_2 \sigma_i \sigma = H_2 \sigma_{h(\sigma)(i)}$ for $i = 1, 2, 3$. (Here if $\tau_1, \tau_2 \in \mathfrak{S}_3$, $(\tau_1 \tau_2)(i) = \tau_2(\tau_1(i))$.) We define a right action of $H_1$ on $\widetilde{V}$ by

$$(x_1 \otimes x_2 \otimes x_3)^\sigma = x^\sigma_{h(\sigma^{-1})(1)} \otimes x^\sigma_{h(\sigma^{-1})(2)} \otimes x^\sigma_{h(\sigma^{-1})(3)}.$$

Let $V = \widetilde{V}^{H_1}$. Then $V$ is a vector space over $k$. We define an action of $G$ on $V$ by

$$(t_1, g)x_1 \otimes x_2 \otimes x_3 = t_1 g^{\sigma_1} x_1 \otimes g^{\sigma_2} x_2 \otimes g^{\sigma_3} x_3.$$

We choose and fix $\sigma \in H_1$ so that $h(\sigma) = (1, 2, 3)$. Then we use $\sigma_1 = 1$, $\sigma_2 = \sigma$, $\sigma_3 = \sigma^2$ as three different imbeddings. If $k_2$ is an $\mathfrak{S}_3$–extension of $k$, this means that we are restricting elements of $\mathrm{Gal}(k_2/k)$ to $k_1$. Also since $\sigma_1 = 1$, the non-trivial element $\nu \in H_2$ corresponds to the transposition $(2, 3)$.

Note that if $x_1, x_2, x_3 \in k$, $(x_1 \otimes x_2 \otimes x_3)^\sigma = x_3 \otimes x_1 \otimes x_2$. This means that if $x = \sum_{i,j,k} x_{ijk} e_{ijk}$,

$$x^\sigma = x^\sigma_{111} e_{111} + x^\sigma_{112} e_{211} + x^\sigma_{121} e_{112} + x^\sigma_{122} e_{212}$$
$$+ x^\sigma_{211} e_{121} + x^\sigma_{212} e_{221} + x^\sigma_{221} e_{122} + x^\sigma_{222} e_{222}.$$

Also
$$x^\nu = x^\nu_{111} e_{111} + x^\nu_{112} e_{121} + x^\nu_{121} e_{112} + x^\nu_{211} e_{211}$$
$$+ x^\nu_{122} e_{122} + x^\nu_{212} e_{221} + x^\nu_{221} e_{212} + x^\nu_{222} e_{222}.$$

Therefore, $x \in V$ is equivalent to

$$x_{111}, x_{222} \in k, \ x_{211}, x_{122} \in k_1,$$
$$x_{122}^{\sigma} = x_{212}, \ x_{122}^{\sigma^2} = x_{221}, \ x_{211}^{\sigma} = x_{121}, \ x_{211}^{\sigma^2} = x_{112}.$$

So $x$ is determined by $x_{111}, x_{222} \in k$ and $x_{122}, x_{211} \in k_1$.

We proved in [1] that $(G, V)$ is a prehomogeneous vector space.

Let $T_2 \subset G_2$ be the subgroup of diagonal matrices, $N_2 \subset G_2$ the subgroup of lower triangular matrices with diagonal entries 1, and $B_2 = T_2 N_2$. We define $T = G_1 \times T_2$, $N = \{1\} \times N_2$. Then $B = TN$ is a Borel subgroup of $G$. We define

$$(1.1) \quad a(t_{21}, t_{22}) = \begin{pmatrix} t_{21} & 0 \\ 0 & t_{22} \end{pmatrix}, \ n(u) = \begin{pmatrix} 1 & 0 \\ u & 1 \end{pmatrix}, \ \tau = \begin{pmatrix} 0 & 1 \\ 1 & 0 \end{pmatrix}.$$

where $t_{21}, t_{22} \in \mathbb{A}_{k_1}^{\times}$, $u \in \mathbb{A}_{k_1}$. We use the same notation for elements $(1, a(t_{21}, t_{22}))$ and $(1, n(u))$. Then $T_2$, $N_2$ consist of elements of the form $a(t_{21}, t_{22})$, $n(u)$ respectively.

The action of $a(t_{21}, t_{22})$ is given by

$$a(t_{21}, t_{22})x = \mathrm{N}_{k_1/k}(t_{21}) x_{111} e_{111} + t_{21} t_{22}^{-1} \mathrm{N}_{k_1/k}(t_{22}) x_{122} e_{122}$$
$$+ t_{22} t_{21}^{-1} \mathrm{N}_{k_1/k}(t_{21}) x_{211} e_{211} + \cdots + \mathrm{N}_{k_1/k}(t_{22}) x_{222} e_{222}.$$

For $t_1 \in \mathrm{GL}(1)_k$, the action is given by the usual multiplication of $t_1$.

By easy computations,

$$n(u) e_{111} = e_{111} + u e_{211} + u^{\sigma} e_{121} + u^{\sigma^2} e_{112}$$
$$+ u^{\sigma} u^{\sigma^2} e_{122} + u u^{\sigma^2} e_{212} + u u^{\sigma} e_{221} + \mathrm{N}_{k_1/k}(u) e_{222},$$
$$n(u) e_{112} = e_{112} + u e_{212} + u^{\sigma} e_{122} + u u^{\sigma} e_{222},$$
$$n(u) e_{121} = e_{121} + u e_{221} + u^{\sigma^2} e_{122} + u u^{\sigma^2} e_{222},$$
$$n(u) e_{122} = e_{122} + u e_{222},$$
$$n(u) e_{211} = e_{211} + u^{\sigma} e_{221} + u^{\sigma^2} e_{212} + u^{\sigma} u^{\sigma^2} e_{222},$$
$$n(u) e_{212} = e_{212} + u^{\sigma} e_{222},$$
$$n(u) e_{221} = e_{221} + u^{\sigma^2} e_{222},$$
$$n(u) e_{222} = e_{222}.$$

Therefore, if $n(u)x = \sum_{i,j,k} y_{ijk} e_{ijk}$,

$$y_{111} = x_{111},$$
$$y_{211} = x_{211} + u x_{111},$$
$$y_{122} = x_{122} + u^{\sigma} u^{\sigma^2} x_{111} + u^{\sigma} x_{211}^{\sigma^2} + u^{\sigma^2} x_{211}^{\sigma},$$
$$y_{222} = x_{222} + \mathrm{N}_{k_1/k}(u) x_{111} + \mathrm{Tr}_{k_1/k}(u^{\sigma} u^{\sigma^2} x_{211}) + \mathrm{Tr}_{k_1/k}(u x_{122}).$$

The element $\tau x$ is obtained by exchanging 1 and 2 in the indices of $e_{ijk}$'s (for example $\tau e_{122} = e_{211}$).

The relative invariant of $(G, V)$ can be constructed in the following manner. For $x = (x_{ijk}) \in \widetilde{V}$, we associate a $2 \times 2$ matrix with entries in the space of linear forms in two variables $v = (v_1, v_2)$ as

$$x \to M_x(v) = v_1 \begin{pmatrix} x_{111} & x_{121} \\ x_{211} & x_{221} \end{pmatrix} + v_2 \begin{pmatrix} x_{112} & x_{122} \\ x_{212} & x_{222} \end{pmatrix}.$$

Then $F_x(v) = \det M_x(v)$ is a quadratic form in $v = (v_1, v_2)$.

By an easy computation,

$$\begin{aligned} F_x(v) = &\, (x_{111} x_{122}^{\sigma^2} - x_{211}^\sigma x_{211}) v_1^2 \\ &+ (x_{111} x_{222} + x_{211}^{\sigma^2} x_{122}^{\sigma^2} - x_{211}^\sigma x_{122}^\sigma - x_{211} x_{122}) v_1 v_2 \\ &+ (x_{211}^{\sigma^2} x_{222} - x_{122}^\sigma x_{122}) v_2^2. \end{aligned}$$

Let $\Delta(x)$ be the discriminant of $F_x(v)$.

We showed in [6] that

$$\Delta((g_1, g_2, g_3)x) = (\det g_1 \det g_2 \det g_3)^2 \Delta(x)$$

for $g_1, g_2, g_3 \in \mathrm{GL}(2)_{k_2}$, $x \in \widetilde{V}$. So if we put $\chi(g) = \mathrm{N}_{k_1/k}(\det g)$ for $g \in \mathrm{GL}(2)_{k_1}$, $\chi$ is a $k$–rational character of $G$ and $\Delta(\phi(g)x) = \chi(g)^2 \Delta(x)$. We proved in [1] that $\Delta(x) \in k[V]$. We define $V_k^{\mathrm{ss}} = \{x \in V_k \mid \Delta(x) \neq 0\}$ and call it the set of semi-stable points. If $x \in V_k \setminus \{0\}$ is not semi-stable, we say it is unstable.

For the rest of this section, we will consider a stratification of $V_k$. Let

$$Y_1 = \{x \in V \mid x_{111} = x_{211} = 0\}, \ Y_1^{\mathrm{ss}} = \{x \in Y_1 \mid x_{122} \neq 0\},$$
$$Z_1 = \{x \in V \mid x_{111} = x_{211} = x_{222} = 0\}, \ Z_1^{\mathrm{ss}} = \{x \in Z_1 \mid x_{122} \neq 0\},$$
$$Y_2 = Z_2 = \{x \in V \mid x_{111} = x_{211} = x_{122} = 0\}, \ Y_2^{\mathrm{ss}} = Z_2^{\mathrm{ss}} = \{x \in Z_2 \mid x_{222} \neq 0\}.$$

We define $S_i = G Y_i^{\mathrm{ss}}$ for $i = 1, 2$.

**Proposition (1.2)** (1) $V_k \setminus \{0\} = V_k^{\mathrm{ss}} \coprod S_{1k} \coprod S_{2k}$.
(2) $S_{ik} \cong G_k \times_{B_k} Y_{ik}^{\mathrm{ss}}$ for $i = 1, 2$.

*Proof.* Let $x \in V_k \setminus \{0\}$ be an unstable point. We first prove that there exists $g \in G_k$ such that $F_{gx}(v) = v_2^2$.

If

$$x_{111} x_{122}^{\sigma^2} - x_{211}^\sigma x_{211} = x_{211}^{\sigma^2} x_{222} - x_{122}^\sigma x_{122} = 0,$$

$F_{gx}(v)$ is a scalar multiple of $v_1 v_2$ and $x \in V_k^{\mathrm{ss}}$. Applying $\tau$ if necessary, we may assume that $x_{111} x_{122}^{\sigma^2} - x_{211}^\sigma x_{211} \neq 0$.

Let $u \in k_1$ and $y = n(u)x$. By an easy computation,

$$F_y(v) = (x_{111} x_{122}^{\sigma^2} - x_{211}^\sigma x_{211}) v_1^2 + (A(x) + 2B(x)u)^{\sigma^2} v_1 v_2 + C(x) v_2^2,$$

where
$$A(x) = x_{111}x_{222} + 2x_{211}x_{122} - \text{Tr}_{k_1/k}(x_{211}x_{122}),$$
$$B(x) = x_{111}x_{122} - x_{211}^{\sigma} x_{211}^{\sigma^2},$$
and we don't have to know what $C(x)$ is.

If $x_{211} \neq 0$, $x_{211}^{\sigma} x_{211}^{\sigma^2} = \text{N}_{k_1/k}(x_{211})x_{211}^{-1} \in k_1^{\times}$. Also $x_{111}, x_{222} \in k \subset k_1$ and $x_{211}, x_{122} \in k_1$. Therefore, we can choose $u \in k_1$ so that the second term becomes zero. Since $x$ is unstable, this means that the third term is zero also. Then for a suitable $t_1 \in k^{\times}$, $F_{t_1 \tau y}(v) = v_2^2$.

Next we prove that for any unstable point $x \in V_k$, there exists $g \in G_k$ such that $gx \in Y_{1k}$. If $x \neq 0$, this means either $gx \in Y_{1k}^{\text{ss}}$ or $Y_{2k}^{\text{ss}}$. We may assume that $F_x(v) = v_2^2$. So the matrix $\begin{pmatrix} x_{111} & x_{211}^{\sigma} \\ x_{211} & x_{122}^{\sigma^2} \end{pmatrix}$ is not invertible. If $x_{111} = 0$, $x_{211} = 0$ also. So $x \in Y_1$. Assume that $x_{111} \neq 0$. Let $u \in k_1$ and $y = n(u)x$. Then $y_{111} = x_{111}$, $y_{211} = x_{211} + ux_{111}$ and the matrix $\begin{pmatrix} y_{111} & y_{211}^{\sigma} \\ y_{211} & y_{122}^{\sigma^2} \end{pmatrix}$ is not invertible either. We can choose $u$ so that $y_{211} = 0$. This implies $y_{122} = 0$. If $y_{222} \neq 0$, $F_y(v)$ is a scalar multiple of $v_1 v_2$ which is a contradiction. Therefore, $y_{222} = 0$. Then $\tau y \in Z_2$.

Next we prove that if $x, y \in Y_{1k}^{\text{ss}}$ and $gx = y$ for $g \in G_k$, there exists $g' \in B_k$ such that $g'x = y$. Note that because of the Bruhat decomposition of $\text{GL}(2)_{k_1}$, we have a decomposition $G_k = B_k \coprod B_k \tau B_k$. If $g \in B_k$, there is nothing to prove. If $g \in B_k \tau B_k$, by replacing $x$ if necessary, we may assume that $g = (t_1, a_2(t_{21}, t_{22})n(u)\tau)$ for certain $t_1 \in k^{\times}$, $t_{21}, t_{22} \in k_1^{\times}$, $u \in k_1$. If $x' = (x'_{ikj}) = \tau x$, $x'_{211} \neq 0$. If $x'_{111} \neq 0$, $y_{111} \neq 0$ which is a contradiction. If $x'_{111} = 0$, since $x'_{211} \neq 0$, $y_{211} \neq 0$ which is a contradiction also. Therefore, we can take the above $g$ from $B_k$.

The argument is similar for $S_2$, and this proves the proposition.

Q.E.D.

## §2 The zeta function and the smoothed Eisenstein series

In this section, we define the zeta function and discuss basic properties of the smoothed Eisenstein series.

Let

$G_{2\mathbb{A}}^0 = \{g \in G_{2\mathbb{A}} \mid |\det g|_1 = 1\}$, $G_{\mathbb{A}}^0 = \mathbb{A}^1 \times G_{2\mathbb{A}}^0$,

$B_{2\mathbb{A}}^1 = \{a(\underline{\mu}_{k_1} t_{21}, \underline{\mu}_{k_1}^{-1} t_{22})n(u) \mid \mu \in \mathbb{R}_+, t_{21}, t_{22} \in \mathbb{A}_{k_1}^1, u \in \mathbb{A}_{k_1}\}$, $B_{\mathbb{A}}^1 = \mathbb{A}^1 \times B_{2\mathbb{A}}^1$.

We define $\widetilde{G}_{\mathbb{A}_k} = \mathbb{R}_+ \times G_{\mathbb{A}_k}^0$. The group $\widetilde{G}_{\mathbb{A}}$ acts on $V_{\mathbb{A}}$ by assuming that $\lambda \in \mathbb{R}_+$ acts by multiplication by $\underline{\lambda}$. Throughout this paper, we express elements $\widetilde{g} \in \widetilde{G}_{\mathbb{A}}$, $g^0 \in G_{\mathbb{A}}^0$ as $\widetilde{g} = (\lambda, t_1, g_2)$, $g^0 = (t_1, g_2)$ where $\lambda \in \mathbb{R}_+$, $t_1 \in \mathbb{A}^1$, and $g_2 \in G_{2\mathbb{A}}^0$. We identify the element $g^0$ with $(1, g^0)$ and $g_2 \in G_{2\mathbb{A}}^0$ with $(1, 1, g_2)$. Let $K_1 \subset \mathbb{A}$ be the unique maximal compact subgroup, and $K_2 \subset G_{2\mathbb{A}}$ the standard maximal compact subgroup. Then $K = K_1 \times K_2$ is a maximal compact subgroup of $G_{\mathbb{A}}^0$. Let $dk$ be the Haar measure on $K$ such that $\int_K dk = 1$. Throughout this paper, we express elements of $B_\mathbb{A}^1$ as

$$b^1 = (t_1, a(\underline{\mu}_{k_1} t_{21}, \underline{\mu}_{k_1}^{-1} t_{22}) n(u)),$$

where $t_1 \in \mathbb{A}^1$, $t_{21}, t_{22} \in \mathbb{A}_{k_1}^1$, $\mu \in \mathbb{R}_+$, $u \in k_1$.

We use the Haar measure $d^\times t'$ on $t' \in \mathbb{A}^1$ or $\mathbb{A}_{k_1}^1$ such that the volume of $\mathbb{A}^1/k^\times$ or $\mathbb{A}_{k_1}^1/k_1$ is 1. If $t = \underline{\lambda} t'$ or $\underline{\lambda}_{k_1} t'$ where $t' \in \mathbb{A}^1$ or $\mathbb{A}_{k_1}^1$, we use $d^\times t = d^\times \lambda d^\times t'$ as the Haar measure on $\mathbb{A}^1$ or $\mathbb{A}_{k_1}^1$. We do not use the subscript $k_1$ for the measure on $\mathbb{A}_{k_1}^1$ because the situation will be clear from the context. We use the measure

$$db^1 = \mu^{-2} d^\times t_1 d^\times \mu d^\times t_{21} d^\times t_{22} du$$

on $B_\mathbb{A}^1$. We use $dg^0 = dk db^1$ as the Haar measure on $G_\mathbb{A}^0$. Let $dg_2$ be the Haar measure on $G_\mathbb{A}^0$ which we defined in §1.1 of [8]. Then $dg^0 = d^\times t_1 dg_2$. We define $d\widetilde{g} = d^\times \lambda dg^0$.

For $\eta > 0$, we define

$$T_{2\eta+}^0 = \{a(\underline{\mu}_{k_1}, \underline{\mu}_{k_1}^{-1}) \mid \mu \in \mathbb{R}_+,\ \mu \geq \eta\}.$$

Let $\Omega \subset B_{2\mathbb{A}}^0$ be a compact subset. We define $\mathfrak{S}^0 = K T_{2\eta+}^0 \Omega$. The set $\mathfrak{S}^0$ is called a Siegel set. It is known that for a suitable choice of $\eta$ and $\Omega$, $\mathfrak{S}^0$ surjects to $G_{2\mathbb{A}}^0/G_{2k}$. Also there exists another compact set $\widehat{\Omega} \subset G_{2\mathbb{A}}^0$ such that $\mathfrak{S}^0 \subset \widehat{\Omega} T_{2\eta+}^0$.

For $x = (x_{ijk})$, $y = (y_{ijk}) \in V$, we define

$$[x,y]' = x_{111} y_{111} + \mathrm{Tr}_{k_1/k}(x_{211} y_{211} + x_{122} y_{122}) + x_{222} y_{222}.$$

This is a non-degenerate bilinear form on $V$ defined over $k$. By an easy computation, $[gx, y] = [x, {}^t gy]$ for $g = (t_1, a_2(t_{21}, t_{22}))$, $\tau$, $n(u)$ and $x, y \in V_\mathbb{A}$. So this is true for all $g \in G$ and $x, y \in V_\mathbb{A}$. We define $[x, y] = [x, \tau y]$. For $\widetilde{g} \in \widetilde{G}$, we define $\widetilde{g}^\iota = (\lambda^{-1}, t^{-1}, \tau {}^t g_2^{-1} \tau)$. This is an involution and the above bilinear form satisfies $[\widetilde{g} x, y] = [x, (\widetilde{g}^\iota)^{-1} y]$ for all $x, y \in V_\mathbb{A}$.

For $\Phi \in \mathscr{S}(V_\mathbb{A})$, we define its Fourier transform by

$$\widehat{\Phi}(x) = \int_{V_\mathbb{A}} \Phi(y) <[x,y]> dy.$$

It is easy to see that the Fourier transform of $\Phi(\widetilde{g} \cdot\ )$ is $\lambda^{-8} \widehat{\Phi}(\widetilde{g}^\iota \cdot\ )$.

**Definition (2.1)** *For any $G_k$-invariant subset $S \subset V_k$ and $\Phi \in \mathscr{S}(V_\mathbb{A})$, we define*

$$\Theta_S(\Phi, \widetilde{g}) = \sum_{x \in S} \Phi(\widetilde{g} x).$$

Note that by Proposition (1.2), if $x \in V_k^{\mathrm{ss}}$, $(x_{111}, x_{211}) \neq 0$ and $(x_{122}, x_{222}) \neq 0$. Therefore, by Lemma (1.2.3) and Lemma (1.2.6) in [8], we get the following lemma.

**Lemma (2.2)** (1) *There exists a slowly increasing function $h(\lambda, \mu)$ such that for any $N_1, N_2 \geq 1$,*

$$\Theta_{V_k^{\mathrm{ss}}}(\Phi, \widetilde{g}) \ll \sup((\lambda \mu)^{-N_1}, (\lambda^3 \mu)^{-N_1}) \sup((\lambda \mu^{-1})^{-N_2}, (\lambda^3 \mu^{-1})^{-N_1}) h(\lambda, \mu).$$

(2) *The function $\Theta_{V_k}(\Phi, \widetilde{g})$ is slowly increasing.*

Note that $x_{211}, x_{122} \in k_1$ and $k_1$ is a vector space of dimension 3 over $k$. So when $x_{211} \neq 0$, we are applying Lemma (1.2.6) in [8] with $m_{ij} = 3$ for a suitable $(i, j)$. Also $\underline{\mu}_{k_1} = \underline{\mu}^{\frac{1}{3}}$. So the multiplication of $\underline{\mu}_{k_1}$ can be considered as a multiplication of an element of $\mathbb{A}$. This is why we get the factor $\lambda^3 \mu$ or $\lambda^3 \mu^{-1}$.

Let $\omega_1, \omega_2$ be characters of $\mathbb{A}_k^\times / k^\times$, $\mathbb{A}_{k_1}^\times / k_1^\times$ respectively. We put $\omega = (\omega_1, \omega_2)$. For this $\omega$, we define a character of $\widetilde{G}_{\mathbb{A}_k} / G_k$ by $\omega(\widetilde{g}) = \omega_1(t_1)\omega_2(\det g_2)$. We define $\delta(\omega_1) = 1$ if $\omega_1$ is a trivial character and $\delta(\omega_1) = 0$ otherwise. We define $\delta(\omega_2)$ similarly.

**Definition (2.3)** *For $\Phi \in \mathscr{S}(V_\mathbb{A})$, $\omega$ as above and a complex variable $s$, we define*

(1) $$Z(\Phi, \omega, s) = \int_{\widetilde{G}_{\mathbb{A}_k}/G_k} \lambda^s \omega(\widetilde{g}) \Theta_{V_k^{ss}}(\Phi, \widetilde{g}) d\widetilde{g},$$

(2) $$Z_+(\Phi, \omega, s) = \int_{\substack{\widetilde{G}_{\mathbb{A}_k}/G_k \\ \lambda \geq 1}} \lambda^s \omega(\widetilde{g}) \Theta_{V_k^{ss}}(\Phi, \widetilde{g}) d\widetilde{g}.$$

By Lemma (2.2), the integral (1) converges absolutely and locally uniformly in a certain right half plane and the integral (2) is an entire function. This means that $(G, V)$ is of complete type (see [8, p. 64]) and therefore, by Theorem (0.3.7) in [8] (which is due to Shintani), $Z(\Phi, \omega, s)$ can be continued meromorphically to the entire plane and satisfies a functional equation

$$Z(\Phi, \omega, s) = Z(\widehat{\Phi}, \omega^{-1}, 8 - s).$$

For $\lambda \in \mathbb{R}_+$ and $\Phi \in \mathscr{S}(V_\mathbb{A})$, we define $\Phi_\lambda(x) = \Phi(\underline{\lambda} x)$. Let

$$J(\Phi, g^0) = \sum_{x \in V_k \setminus V_k^{ss}} \widehat{\Phi}((g^0)^\iota x) - \sum_{x \in V_k \setminus V_k^{ss}} \Phi(g^0 x),$$

$$I^0(\Phi, \omega) = \int_{G_\mathbb{A}^0 / G_k} \omega(g^0) J(\Phi, g^0) dg^0,$$

$$I(\Phi, \omega, s) = \int_0^1 \lambda^s I^0(\Phi_\lambda, \omega) d^\times \lambda.$$

Then by the Poisson summation formula,

$$Z(\Phi, \omega, s) = Z_+(\Phi, \omega, s) + Z_+(\widehat{\Phi}, \omega^{-1}, 8 - s) + I(\Phi, \omega, s).$$

We study the last term in §3. For that purpose, we briefly discuss basic properties of the smoothed Eisenstein series. For $z = (z_1, z_2) \in \mathbb{C}^2$ such that $z_1 + z_2 = 0$, we define $a(t_{21}, t_{22})^z = |t_{21}|_{k_1}^{z_1} |t_{22}|_{k_1}^{z_2}$. Let $\rho = (\frac{1}{2}, -\frac{1}{2})$. We can consider $\rho$ as half the sum of positive weights. Let $\psi(z)$ be an entire function which is rapidly decreasing with respect to $\text{Im}(z)$ on any vertical strip. Moreover, we assume that

$$\psi(z_1, z_2) = \psi(z_2, z_1), \quad \psi(\rho) = \psi\left(\frac{1}{2}, -\frac{1}{2}\right) \neq 0.$$

Let $\mathscr{E}(g_2, w, \psi)$ be the smoothed Eisenstein series defined for $g_2 \in G_{2\mathbb{A}}^0$ and a complex variable $w$. For the definition of $\mathscr{E}(g_2, w, \psi)$, see [3, p. 172], [5, p. 526], [7, p. 368], or [8, p. 75], (we are applying the definition in [5], [7], [8] replacing $k$ by $k_1$). Since the dependency on $\psi$ is not so crucial, we drop $\psi$ and use the notation $\mathscr{E}(g_2, w)$ instead.

Let $s = z_1 - z_2$. Since $z_1 + z_2 = 0$, $\psi$ can be considered as a function of $s$. With this understanding, $\rho$ corresponds to $s = 1$. We define

$$(2.4) \qquad \Lambda(w; s) = \frac{\psi(s)}{w - s}.$$

We recall properties of $\mathscr{E}(g_2, w)$. Let

$$\phi_{k_1}(s) = Z_{k_1}(s) Z_{k_1}(s+1)^{-1},$$
$$\mathfrak{R}_{k_1} = \operatorname*{Res}_{s=1} Z_{k_1}(s),$$
$$\varrho_{k_1} = \operatorname*{Res}_{s=1} \phi_{k_1}(s) = \frac{\mathfrak{R}_{k_1}}{Z_{k_1}(2)},$$
$$\mathfrak{V}_{k_1 2} = \varrho_{k_1}^{-1},$$

where $Z_{k_1}(s)$ is the Dedekind zeta function of the field $k_1$ including the Gamma factors. The function $Z_{k_1}(s)$ is defined in [4, p. 129] or [8, p. xii] (for the field $k$). It is well known that $\mathfrak{V}_{k_1 2}$ is the volume of $G_{2\mathbb{A}}^0 / G_{2k}$.

The following Proposition is well known.

**Proposition (2.5)** (1) If $f(g_2)$ is a function of $g_2 \in G_{2\mathbb{A}}^0 / G_{2k}$ such that there is a constant $r < 2$ and

$$f(k_2 a(\underline{\mu}_{k_1}, \underline{\mu}_{k_1}^{-1}) n(u)) \ll \mu^r$$

for $k_2 a(\underline{\mu}_{k_1}, \underline{\mu}_{k_1}^{-1}) n(u) \in \mathfrak{S}^0$, the integral

$$\int_{G_{2\mathbb{A}}^0 / G_{2k}} f(g_2) \mathscr{E}(g_2, w) dg_2$$

becomes a holomorphic function for $\operatorname{Re}(w) \geq 1 - \epsilon$ for a constant $\epsilon > 0$ except possibly for a simple pole at $w = 1$ with residue

$$\varrho_{k_1} \int_{G_{2\mathbb{A}}^0 / G_{2k}} f(g_2) dg_2.$$

(2) If $f(g_2)$ is a slowly increasing function of $g_2 \in G_{2\mathbb{A}}^0 / G_{2k}$, the integral

$$\int_{G_{2\mathbb{A}}^0 / G_{2k}} f(g_2) \mathscr{E}(g_2, w) dg_2$$

becomes a holomorphic function on a certain right half plane.

(3) If $\omega_2$ is a character of $\mathbb{A}_{k_1}^\times / k_1^\times$,

$$\int_{G_{2\mathbb{A}}^0 / G_{2k}} \omega_2(\det g_2) \mathscr{E}(g_2, w) dg_2 = \delta(\omega_2) \Lambda(w; \rho).$$

Let $\mathscr{E}_0(g_2, w)$, $\widetilde{\mathscr{E}}(g_2, w)$ be the constant term and the non-constant term of $\mathscr{E}(g_2, w)$ respectively. It is well known that

$$(2.6) \qquad \mathscr{E}_0(b^1, w) = \frac{1}{2\pi\sqrt{-1}} \int_{\mathrm{Re}(s)=r>1} (\mu^{1+s} + \mu^{1-s}\phi_{k_1}(s))ds.$$

Also on any vertical strip in $\mathrm{Re}(w) > 0$, $\widetilde{\mathscr{E}}(b^1, w)$ is holomorphic and

$$(2.7) \qquad \widetilde{\mathscr{E}}(b^1, w) \ll \mu^{1-2l}.$$

For the proof of Proposition (2.5), see [3, p. 172], [5, p. 527], or §3.4 of [8].

If $f(w), g(w)$ are meromorphic functions of $w$, we use the notation $f \sim g$ if $f - g$ can be continued meromorphically to $\mathrm{Re}(w) > 1 - \epsilon$ for a constant $\epsilon > 0$ and is holomorphic at $w = 1$. We define

$$I(\Phi, \omega, w) = \int_{G_\mathbb{A}^0/G_k} \omega(g^0)\mathscr{E}(g_2, w)J(\Phi, g^0)dg^0.$$

By Lemma (2.2) and Proposition (2.5),

$$(2.8) \qquad I(\Phi, \omega, w) \sim \varrho_{k_1}\Lambda(w; \rho)I^0(\Phi, \omega).$$

§3 The principal part formula

In this section, we determine the principal part of the zeta function. For that purpose, we first have to define distributions which appear in the inductive computation.

**Definition (3.1)** (1) *Let $\Psi_1, \Psi_2$ be Schwartz–Bruhat functions on $\mathbb{A}$, $\mathbb{A}_{k_1}$ respectively. Then for $t_1 \in \mathbb{A}^\times$, $t_2 \in k_{k_1}^\times$, we define*

$$\Theta_1(\Psi_1, t_1) = \sum_{x \in k^\times} \Psi_1(t_1 x),$$

$$\Theta_{k_1 1}(\Psi_2, t_2) = \sum_{x \in k_1^\times} \Psi_2(t_2 x).$$

(2) *Let $\omega_1, \omega_2$ be characters of $\mathbb{A}^\times/k^\times$, $\mathbb{A}_{k_1}^\times/k_1^\times$ respectively. For a complex variable $s$, we define*

$$\Sigma_1(\Psi_1, \omega_1, s) = \int_{\mathbb{A}^\times/k^\times} |t_1|^s \omega_1(t_1)\Theta_1(\Psi_1, t_1)d^\times t_1,$$

$$\Sigma_{k_1 1}(\Psi_2, \omega_2, s) = \int_{\mathbb{A}_{k_1}^\times/k_1^\times} |t_1|_{k_1}^s \omega_2(t_2)\Theta_{k_1 1}(\Psi_2, t_2)d^\times t_2.$$

These are more or less the Dedekind zeta functions of $k$, $k_1$ respectively. They can be continued meromorphically to the entire plane with possible simple poles at $s = 0, 1$.

Let

$$\Sigma_{k_1 1}(\Psi_2, \omega_2, s) = \sum_{i=-1}^{\infty} \Sigma_{k_1 1,(i)}(\Psi_2, \omega_2, s_0)(s - s_0)^i \tag{3.2}$$

be the Laurent expansion at $s = s_0$.

Let $\Phi \in \mathscr{S}(V_\mathbb{A})$. We define a Schwartz–Bruhat function $R_1\Phi$ on $Z_{1\mathbb{A}} \cong \mathbb{A}_{k_1}$ by

$$R_1\Phi(x_{122}) = \int_\mathbb{A} \Phi(0, 0, x_{122}, x_{222}) dx_{222}. \tag{3.3}$$

We define a Schwartz–Bruhat function $R_2\Phi$ on $Z_{2\mathbb{A}} \cong \mathbb{A}$ by just restricting $\Phi$ to $Z_{2\mathbb{A}}$.

**Definition (3.4)** (1) $\delta_\#(\omega) = \delta(\omega_1)\delta(\omega_2)$.
(2) $\delta_1(\omega) = 1$ if $\omega_1(q_1)\omega_2(q_1)^{-1} = \omega_2(q_2^2 N_{k_1/k}(q_2)^{-1}) = 1$ for all $q_1 \in \mathbb{A}^\times$, $q_2 \in \mathbb{A}_{k_1}^\times$ and $\delta_1(\omega) = 0$ otherwise
(3) $\delta_2(\omega) = 1$ if $\omega_1(N_{k_1/k}(q)) = \omega_2(q) = 1$ for all $q \in \mathbb{A}_{k_1}^\times$ and $\delta_2(\omega) = 0$ otherwise

As in [5] or [8], for $\Phi \in \mathscr{S}(V_\mathbb{A})$ and $\omega = (\omega_1, \omega_2)$ as above, we define

$$M_\omega \Phi(x) = \int_K \omega(k)\Phi(kx) dk. \tag{3.5}$$

Note that $M_\omega M_\omega \Phi = M_\omega \Phi$ and $Z(\Phi, \omega, s) = Z(M_\omega \Phi, \omega, s)$. So in order to study $Z(\Phi, \omega, s)$, we may assume that $M_\omega \Phi = \Phi$. The following theorem is the main theorem of this paper.

**Theorem (3.6)** *Suppose $M_\omega \Phi = \Phi$. Then*

$$Z(\Phi, \omega, s) = Z_+(\Phi, \omega, s) + Z_+(\widehat{\Phi}, \omega^{-1}, 8 - s)$$
$$+ \mathfrak{V}_{k_1 2}\delta_\#(\omega)\left(\frac{\widehat{\Phi}(0)}{s - 8} - \frac{\Phi(0)}{s}\right)$$
$$+ \delta_2(\omega)\left(\frac{\Sigma_1(R_2\widehat{\Phi}, \omega_1^{-1}, 2)}{s - 6} - \frac{\Sigma_1(R_2\Phi, \omega_1, 2)}{s - 2}\right)$$
$$+ \delta_1(\omega)\left(-\frac{3\Sigma_{k_1 1,(-1)}(R_1\widehat{\Phi}, \omega_2^{-1}, 1)}{(s - 4)^2} + \frac{\Sigma_{k_1 1,(0)}(R_1\widehat{\Phi}, \omega_2^{-1}, 1)}{s - 4}\right)$$
$$- \delta_1(\omega)\left(\frac{3\Sigma_{k_1 1,(-1)}(R_1\Phi, \omega_2, 1)}{(s - 4)^2} + \frac{\Sigma_{k_1 1,(0)}(R_1\Phi, \omega_2, 1)}{s - 4}\right).$$

*Proof.* We devote the rest of this paper to the proof of this theorem. We use Proposition (2.5) and (2.8) to separate contributions from unstable points. The contribution from each stratum is given by the following integral

$$I_i(\Phi, \omega, w) = \int_{G_\mathbb{A}^0/G_k} \omega(g^0)\Theta_{S_{ik}}(\Phi, g^0)\mathscr{E}(g_2, w) dg^0 \text{ for } i = 1, 2. \tag{3.7}$$

Since $\mathscr{E}(g_2, w) = \mathscr{E}(g_2^l, w)$, by Proposition (2.5),

(3.8) $$I(\Phi, \omega, w) = \delta_{\#}(\omega)\Lambda(w; \rho)(\widehat{\Phi}(0) - \Phi(0))$$
$$+ \sum_{i=1,2}(I_i(\widehat{\Phi}, \omega^{-1}, w) - I_i(\Phi, \omega, w)).$$

We first consider $I_1(\Phi, \omega, w)$.

**Proposition (3.9)**

$$I_1(\Phi, \omega, w) \sim \frac{\delta_1(\omega)}{2\pi\sqrt{-1}} \int_{\mathrm{Re}(s)=r>1} \Sigma_{k_1 1}(R_1\Phi, \omega_2, s)\phi_{k_1}(s)\Lambda(w;s)ds.$$

*Proof.*
$$I_1(\Phi, \omega, w) = \int_{G_{\mathbb{A}}^0/G_k} \omega(g^0)\Theta_{S_{1k}}(\Phi, g^0)\mathscr{E}(g_2, w)dg^0$$
$$= \int_{G_{\mathbb{A}}^1/B_k} \omega(g^0)\Theta_{Y_{1k}^{\mathrm{ss}}}(\Phi, g^0)\mathscr{E}(g_2, w)dg^0$$
$$= \int_{B_{\mathbb{A}}^1/B_k} \omega(b^1)\Theta_{Y_{1k}^{\mathrm{ss}}}(\Phi, b^1)\mathscr{E}(b^1, w)db^1.$$

**Lemma (3.10)**
$$\int_{B_{\mathbb{A}}^1/B_k} \omega(b^1)\Theta_{Y_{1k}^{\mathrm{ss}}}(\Phi, b^1)\widetilde{\mathscr{E}}(b^1, w)db^1 \sim 0.$$

*Proof.* By Lemma (1.2.8) in [8], for any $N \geq 1$,
$$\Theta_{Y_{1k}^{\mathrm{ss}}}(\Phi, b^1) \ll \mu^N \sup(1, \mu).$$

So on any vertical strip in $\mathrm{Re}(w) > 0$, for any $N \geq 1$, $l > 2$,
$$\Theta_{Y_{1k}^{\mathrm{ss}}}(\Phi, b^1)\widetilde{\mathscr{E}}(b^1, w) \ll \mu^{N+1-2l}\sup(1,\mu).$$

This proves the lemma.
$$\text{Q.E.D.}$$

Note that $\mathscr{E}_0(b^1, w)$ depends only on $\mu$. By substituting $q = t_1 t_{21} t_{22}^{-1} \mathrm{N}_{k_1/k}(t_{22})$, we get

$$\int_{B_{\mathbb{A}}^1/B_k} \omega(b^1)\Theta_{Y_{1k}^{\mathrm{ss}}}(\Phi, b^1)\mathscr{E}_0(b^1, w)db^1$$
$$= \delta_1(\omega)\int_{\mathbb{R}_+ \times \mathbb{A}_{k_1}^{\times}/k_1^{\times}} \omega_2(q)\Theta_{k_1 1}(R_1\Phi, \underline{\mu}_{k_1}^{-1}q)\mathscr{E}_0(a(\underline{\mu}_{k_1}, \underline{\mu}_{k_1}^{-1}), w)\mu^{-1}d^{\times}\mu d^{\times}q.$$

It is easy to see that

$$\int_{\mathbb{R}_+ \times \mathbb{A}_{k_1}^{\times}/k_1^{\times}} \omega_2(q)\Theta_{k_1 1}(R_1\Phi, \underline{\mu}_{k_1}^{-1}q)\mu^{\pm s}d^{\times}\mu d^{\times}q$$
$$= \Sigma_{k_1 1}(R_1\Phi, \omega_2, \mp s)$$

if Re($s$) < −1, Re($s$) > 1 respectively. Since

$$\frac{1}{2\pi\sqrt{-1}} \int_{\text{Re}(s)=r<-1} \Sigma_{k_1 1}(R_1\Phi, \omega_2, -s)\Lambda(w; s)ds \sim 0,$$

we get the proposition.

Q.E.D.

**Proposition (3.11)**

$$I_2(\Phi, \omega, w) \sim \varrho_{k_1}\delta_\#(\omega)\Lambda(w; \rho)\Sigma_1(R_2\Phi, \omega_1, 2).$$

*Proof.*

$$I_2(\Phi, \omega, w) = \int_{G_\mathbb{A}^0/G_k} \omega(g^0)\Theta_{S_{2k}}(\Phi, g^0)\mathscr{E}(g_2, w)dg^0$$

$$= \int_{G_\mathbb{A}^1/B_k} \omega(g^0)\Theta_{Z_{2k}^{ss}}(\Phi, g^0)\mathscr{E}(g_2, w)dg^0$$

$$= \int_{B_\mathbb{A}^1/B_k} \omega(b^1)\Theta_{Z_{2k}^{ss}}(\Phi, b^1)\mathscr{E}(b^1, w)db^1$$

$$= \int_{B_\mathbb{A}^1/B_k} \omega(b^1)\Theta_{Z_{2k}^{ss}}(\Phi, b^1)\mathscr{E}_0(b^1, w)db^1.$$

The last step is because $\Theta_{Z_{2k}^{ss}}(\Phi, b^1)$ does not depend on $u$.

By substituting $q = t_1 N_{k_1/k}(t_{22})$, we get

$$\int_{B_\mathbb{A}^1/B_k} \omega(b^1)\Theta_{Z_{2k}^{ss}}(\Phi, b^1)\mathscr{E}_0(b^1, w)db^1$$

$$= \delta_2(\omega) \int_{\mathbb{R}_+ \times \mathbb{A}_{k_1}^\times / k_1^\times} \omega_1(q)\Theta_1(R_2\Phi, \underline{\mu}_{k_1}^{-3}q)\mathscr{E}_0(a(\underline{\mu}_{k_1}, \underline{\mu}_{k_1}^{-1}), w)\mu^{-2}d^\times\mu d^\times q.$$

Since $\underline{\mu}_{k_1}^{-3} = \underline{\mu}^{-1}$, It is easy to see that

$$\int_{\mathbb{R}_+ \times \mathbb{A}_{k_1}^\times / k_1^\times} \omega_1(q)\Theta_1(R_2\Phi, \underline{\mu}_{k_1}^{-3}q)\mu^{1\pm s}d^\times\mu d^\times q$$

$$= \Sigma_1(R_2\Phi, \omega_1, 1 \mp s)$$

if Re($s$) < −1, Re($s$) > 1 respectively. Therefore,

$$\delta_2(\omega) \int_{\mathbb{R}_+ \times \mathbb{A}_{k_1}^\times / k_1^\times} \omega_1(q)\Theta_1(R_2\Phi, \underline{\mu}_{k_1}^{-3}q)\mathscr{E}_0(a(\underline{\mu}_{k_1}, \underline{\mu}_{k_1}^{-1}), w)\mu^{-2}d^\times\mu d^\times q$$

$$\sim \frac{\delta_2(\omega)}{2\pi\sqrt{-1}} \int_{\text{Re}(s)=r>1} \Sigma_1(R_2\Phi, \omega_1, 1+s)\phi_{k_1}(s)\Lambda(w; s)ds$$

$$\sim \varrho_{k_1}\delta_2(\omega)\Lambda(w; \rho)\Sigma_1(R_2\Phi, \omega_1, 2).$$

Q.E.D.

Let
$$J_1(\Phi,\omega) = \delta_\#(\omega)(\widehat{\Phi}(0) - \Phi(0)) \tag{3.12}$$
$$+ \varrho_{k_1}\delta_2(\omega)(\Sigma_1(R_2\widehat{\Phi},\omega_1^{-1},2) - \Sigma_1(R_2\Phi,\omega_1^{-1},2)),$$
$$J_2(\Phi,\omega,s) = \delta_1(\omega)(\Sigma_{k_11}(R_1\widehat{\Phi},\omega_2^{-1},s) - \Sigma_{k_11}(R_1\Phi,\omega_2,s)).$$

Then by (3.8) and Propositions (3.9), (3.11),
$$I(\Phi,\omega,w) \sim J_1(\Phi,\omega)\Lambda(w;\rho) + \frac{1}{2\pi\sqrt{-1}}\int_{\mathrm{Re}(s)=r>1} J_2(\Phi,\omega,s)\phi_{k_1}(s)\Lambda(w;s)ds.$$

Therefore, by Wright's principle (see [7, p. 373] or §3.7 of [8] replacing $k$ by $k_1$), $J_2(\Phi,\omega,s)$ is holomorphic at $s=1$ and
$$I(\Phi,\omega,w) \sim J_1(\Phi,\omega)\Lambda(w;\rho) + \varrho_{k_1}J_2(\Phi,\omega,1)\Lambda(w;\rho).$$

This implies
$$\delta_1(\omega)\Sigma_{k_11,(-1)}(R_1\widehat{\Phi},\omega_2^{-1},1) = \delta_1(\omega)\Sigma_{k_11,(-1)}(R_1\Phi,\omega_2,1),$$
$$J_2(\Phi,\omega,1) = \delta_1(\omega)(\Sigma_{k_11,(0)}(R_1\widehat{\Phi},\omega_2^{-1},1) - \Sigma_{k_11,(0)}(R_1\Phi,\omega_2,1)).$$

Therefore,
$$I^0(\Phi,\omega) = \mathfrak{V}_{k_12}\delta_\#(\omega)(\widehat{\Phi}(0) - \Phi(0)) \tag{3.13}$$
$$+ \delta_2(\omega)(\Sigma_1(R_2\widehat{\Phi},\omega_1^{-1},2) - \Sigma_1(R_2\Phi,\omega_1,2))$$
$$+ \delta_1(\omega)(\Sigma_{k_11,(0)}(R_1\widehat{\Phi},\omega_2^{-1},1) - \Sigma_{k_11,(0)}(R_1\Phi,\omega_2,1)).$$

It is easy to see that
$$\Phi_\lambda(0) = \Phi(0),\ \widehat{\Phi_\lambda}(0) = \lambda^{-8}\widehat{\Phi}(0),$$
$$\Sigma_1(R_2\Phi_\lambda,\omega_1,2) = \lambda^{-2}\Sigma_1(R_2\Phi,\omega_1,2),$$
$$\Sigma_1(R_2\widehat{\Phi_\lambda},\omega_1^{-1},2) = \lambda^{-6}\Sigma_1(R_2\widehat{\Phi},\omega_1^{-1},2).$$

Since
$$\Sigma_{k_11}(R_1\Phi_\lambda,\omega_2,s) = \lambda^{-1-3s}\Sigma_{k_11}(R_1\Phi,\omega_2,s),$$
$$\Sigma_{k_11}(R_1\widehat{\Phi_\lambda},\omega_2^{-1},s) = \lambda^{-7+3s}\Sigma_{k_11}(R_1\widehat{\Phi},\omega_2^{-1},s),$$

we get
$$\Sigma_{k_11,(0)}(R_1\Phi_\lambda,\omega_2,1) = \lambda^{-4}\Sigma_{k_11,(0)}(R_1\Phi,\omega_2,1)$$
$$- 3\lambda^{-4}\log\lambda\Sigma_{k_11,(-1)}(R_1\Phi,\omega_2,1),$$
$$\Sigma_{k_11,(0)}(R_1\widehat{\Phi_\lambda},\omega_2^{-1},1) = \lambda^{-4}\Sigma_{k_11,(0)}(R_1\widehat{\Phi},\omega_2^{-1},1)$$
$$+ 3\lambda^{-4}\log\lambda\Sigma_{k_11,(0)}(R_1\widehat{\Phi},\omega_2^{-1},1).$$

Then Theorem (3.6) follows by integrating $\lambda^s I^0(\Phi_\lambda,\omega)$ over $s \in [0,1]$ using (3.13).
$$\text{Q.E.D.}$$

**Remark (3.14)** For a place $v$, let $k_v$ be the corresponding completion of $k$. If $\Phi = \otimes\Phi_v$ where $\Phi_v$ is a Schwartz–Bruhat function on $V_{k_v}$ and $\Phi_v$ is compactly supported in $V_{k_v}^{ss}$ for a certain infinite place $v$, $\Sigma_{k_11,(-1)}(R_2\Phi,\omega_2,1) = 0$. But if $M_\omega\Phi = \Phi$,
$$\Sigma_{k_11,(-1)}(R_2\widehat{\Phi},\omega_2^{-1},1) = \Sigma_{k_11,(-1)}(R_2\Phi,\omega_2,1) = 0$$
also. Therefore, all the poles of the associated Dirichlet series for this prehomogeneous vector space are simple.

Akihiko Yukie
Oklahoma State University
Mathematics Department
401 Math Science
Stillwater OK 74078-1058 USA
yukie@math.okstate.edu